\documentclass[12pt]{article}
\usepackage{amssymb,amsmath, amscd,ifthen}
\usepackage[cp1251]{inputenc}
\usepackage[english]{babel}
\usepackage{indentfirst}
\usepackage[affil-it]{authblk}

\newtheorem{lem}{Lemma}
\newtheorem{thm}{Theorem}

\title{Recency-based preferential attachment models}

\author{Liudmila Ostroumova Prokhorenkova and  Egor Samosvat}

\affil{Yandex, Moscow, Russia}

\date{}

\newcommand{\E}{\mathrm{E}}

\newcommand{\BA}{Barab\'asi--Albert}

\newcommand{\Prob}{\mathrm{P}}

\begin{document}

\newenvironment{Proof}{\noindent{\it Proof.\,}}{\hfill$\Box$}

\maketitle

\begin{abstract}

Preferential attachment models were shown to be very effective in predicting such important properties of real-world networks as the power-law degree distribution, small diameter, etc.
However, they do not allow to model the so-called \textit{recency property}. Recency property reflects the fact that in many real networks vertices tend to connect to other vertices of similar age. This fact motivated us to introduce and analyze a new class of models --- recency-based models. This class is a generalization of fitness models, which were suggested by Bianconi and Barab\'{a}si. Bianconi and Barab\'{a}si extended preferential attachment models with pages' inherent quality or \textit{fitness} of vertices.
To additionally reflect a recency property, it is reasonable to generalize fitness models by adding a recency factor to the attractiveness function.
This means that pages are gaining incoming links according to
their \textit{attractiveness}, which is determined by the incoming
degree of the page (current popularity), its \textit{inherent quality} (some
page-specific constant) and age (new pages are gaining new links more
rapidly).

In this paper, we rigorously analyze the degree distribution in the most realistic recency-based model. Also, we prove that this model does reflect the recency property.

\vspace{20pt}

\end{abstract}

\section{Introduction}

Numerous models have been suggested to reflect and predict the growth
of the Web \cite{Boccaletti06,Bonato04,Kumar00}, the most well-known ones are preferential attachment models.
One of the first attempts to propose a realistic mathematical model of
the Web growth was made in~\cite{Barabasi99}. The main idea is to
start with the assumption that new pages often link
to old popular pages. Barab\'{a}si and Albert defined a graph
construction stochastic process, which is a Markov chain of graphs,
governed by the \emph{preferential attachment}. At each step in
the process, a new vertex is added to the graph and is joined to $m$
different vertices already existing in the graph that are chosen with
probabilities proportional to their incoming degree (the measure of
popularity). This model successfully explained some properties of the
Web graph like its small diameter and power-law
distribution of incoming degrees. Later, many modifications to the \BA~model have been
proposed, e.g., \cite{Buckley04,Cooper03,Holme02}, in order to more
accurately depict these but also other properties (see \cite{Albert02,Bollobas03} for details).

It was noted by Bianconi and Barab\'{a}si in \cite{Bianconi01} that in
real networks some vertices are gaining new incoming links not only because of
their incoming degree (popularity), but also because of their own intrinsic properties.
Motivated by this observation, Bianconi and Barab\'{a}si extended preferential attachment models with pages'
inherent quality or \textit{fitness} of vertices.
When a new vertex is added to the graph, it is joined to some already existing
vertices that are chosen with probabilities proportional to the product
of their fitness and incoming degree.


One of the main drawbacks of these models is that they pay too much attention to old
pages and do not realistically explain how links pointing to newly-created pages appear.
For example, most new media pages like news and blog posts are popular only for a short period of time, i.e., such pages are mostly cited and
visited for several days after they appeared.  In \cite{2012model} a \textit{recency property} was introduced,
which reflects the fact that new media pages tend to connect to other media pages of similar age.
Namely, for the media related part of the Web it was shown that $e(T)$ --- the fraction of edges connecting the pages whose age difference is greater than $T$ ---
decreases exponentially fast.

Although preferential attachment models reflect some important properties of real-world networks,
they do not allow to model the recency property.
Here we discuss recency-based models ---
a generalization of fitness models, where a recency factor is added to the attractiveness function. This means that pages are gaining incoming links according to
their \textit{attractiveness}, which is determined by the incoming
degree of the page (current popularity), its \textit{inherent quality} (some
page-specific constant) and age (new pages are gaining new links more
rapidly).

The idea of adding the recency factor to the attractiveness function was first suggested in~\cite{2012model}.
In this paper, we propose a new formal definition of the model from~\cite{2012model}.
Also, this rigorous definition allows us to theoretically analyze different properties of the recency-based models more thoroughly using a combinatorial approach.
Our analysis shows that for the considered models the power-law distribution of inherent quality leads to the power-law degree distribution.
We also rigorously analyze the recency property, i.e., the behavior of~$e(T)$.

\section{Motivation}\label{Motivation}

In \cite{2012model} a model of the Media Web evolution has been proposed. The idea was to generalize the preferential attachment in the way that the probability to cite a page $p$ is
proportional to the \emph{attractiveness} of $p$, which is some function
of $d(p)$ (current degree of $p$), $q(p)$ (intrinsic quality of $p$), and $a(p)$ (current age of $p$). Different attractiveness functions were considered in~\cite{2012model}:

$$\mathrm{attr}(p) =  q(p)^{\alpha_1} \cdot d(p)^{\alpha_2} \cdot e^{-\frac{a(p)}{\tau}\cdot \alpha_3}\,,$$
where $(\alpha_1, \alpha_2, \alpha_3) \in \{0,1\}^3$ and $\tau$ corresponds to the mean lifetime of the decaying attractiveness.
For example, $\mathrm{attr}(p) = d(p)$ leads to preferential
attachment, while $\mathrm{attr}(p) = q(p) \cdot d(p)$ leads to
fitness model.

To depict the recency property of the Media Web one has to include the recency factor $e^{-\frac{a(p)}{\tau}}$ in the attractiveness  function. Further mean-field approximation analysis and computer simulations~\cite{2012model} showed that in order to have the power-law with a realistic exponent, attractiveness  function $\mathrm{attr}(p) = q(p) \, e^{-\frac{a(p)}{\tau}}$ should be chosen.
Moreover, the distribution of qualities $q$ should follow the power law.
That was also confirmed by the analysis of the likelihood of a real data given the model.

Note that some other recency factors have been previously proposed.
In~\cite{dorogovtsev2002evolution,dorogovtsev2000evolution} attractiveness function $\mathrm{attr}(p) = d(p) \cdot a(p)^{-\alpha}$ was studied using the mean-field approximation and computer simulations.
It was shown that the degree distribution follows the power law for $\alpha< 1$ and is exponential for $\alpha > 1$.
Of course, it is also interesting to analyze $\mathrm{attr}(p) = d(p) \cdot a(p)^{-\alpha}$ in a mathematically rigorous way but this is beyond the scope of this paper.

Thus, in this paper we mostly focus on the attractiveness function $q(p) \cdot e^{-\frac{a(p)}{\tau}}$ and the goal is to formalize this model and then analyze it rigorously.
In addition, our analysis allows to understand under which assumptions the conclusions made in~\cite{2012model} hold.




\section{Model}\label{Model}

In this section, we formalize the model introduced in \cite{2012model}.
We construct a sequence of random graphs $\{G_n\}$.
This sequence has the following parameters: a positive integer constant $m$ (vertex outdegree) and an integer function $N(n)$.
We also need a sequence of mutually independent random variables $\zeta_1, \zeta_2, \ldots$ with some given distribution taking positive values.
Each graph $G_n$ is defined according to its own constructing procedure which is based on the idea of preferential attachment.

Let us now define the random graph $G_n$. At the beginning of the constructing process we have two vertices and one edge between them (graph $\tilde{G}_2^n$). The first two vertices have inherent qualities $q(1):=\zeta_1$ and $q(2):=\zeta_2$. At the $t+1$-th step ($2 \le t \le n-1$) one vertex and $m$ edges are added to $\tilde{G}_t^n$. New vertex $t+1$ has an inherent quality $q(t+1):=\zeta_{t+1}$. New edges are drawn independently and they go from the new vertex to previous vertices. For each edge the probability that it goes to a vertex $i$ ($1 \le i \le t$) is equal~to
$$
\frac{\mathrm{attr}_t(i)}{\sum_{j=1}^t \mathrm{attr}_t(j)},
$$
where
$$
\mathrm{attr}_t(i) = q(i) \, e^{-\frac {t-i}{N(n)}}.
$$
According to the definition, loops are not allowed, although multiple edges may appear.

It is important to note that, in contrast to standard definitions of preferential attachment models, in our case a graph $G_n$ cannot be obtained from a graph $G_{n-1}$. Each graph has its own constructing procedure which is based on preferential attachment.
The reason is that, in contrast to~\cite{2012model}, the mean lifetime of the decaying attractiveness $N(n)$ varies with $n$.
This unusual definition allows us to rigorously analyze both the power-law degree distribution and the behavior of $e(T)$, which is the fraction of edges connecting vertices $i$ and $j$ with $|i-j|>T$.
Further we assume that $N(n) \to \infty$ as $n$ grows.
This allows us to analyze the fraction of vertices of degree $d = d(n)$ when $d(n)$ may grow with $n$.


The choice of the attractiveness function is motivated by the empirical results obtained in  \cite{2012model}.
However, in this paper we also consider the attractiveness function
$$
\mathrm{attr}_t(i) = q(i) \, I[i > t - N(n)]
$$
which approximates $q(i) \, e^{-\frac {t-i}{N(n)}}$.
We do this because
both attractiveness functions  are similar in terms of the degree distribution, but the theoretical analysis of the attractiveness function $ q(i) \, I[i > t - N(n)]$ is less complicated, therefore it can be considered as the natural first step.
In addition, the attractiveness function $q(i) \, I[i > t - N(n)]$ has its own practical intuition, but unfortunately it cannot model the recency property well (see Section~\ref{sec:recency} for the discussion).
In this paper we analyze both attractiveness functions.

Let us remark that according to the definition of the model the attractiveness of vertices decays rapidly with time.
Thus, two vertices with a big age difference are unlikely to be connected and the diameter of a network grows as $n / N(n)$.
In this sense, we always observe a chain structure, as it was noted in~\cite{dorogovtsev2000evolution} for $\alpha>1$.
However, if one considers a subgraph induced by $\sim N(n)$ consecutive vertices,
it will be similar to a standard scale-free network.
Further we omit $n$ in the notation $N(n)$.

\section{Attractiveness function $q(i) I[i > t - N]$}

In this section, we assume that the attractiveness function of a vertex $i$ is $\mathrm{attr}_t(i) = q(i) \, I[i > t - N]$. The indicator function means that a vertex $i$ accumulates incoming edges only during the next $N$ steps after its appearance and we call this period a \textit{lifespan} of a vertex. We say that during this lifespan a vertex is \textit{alive}, after this period a vertex \textit{dies}.

As discussed in Section~\ref{Motivation}, we also assume that the random variables $\zeta_1, \zeta_2, \ldots$ have the Pareto distribution with the density function $f(x) = \frac{\gamma  a^{\gamma} I[x > a]}{x^{\gamma+1}}$, where $\gamma > 1$, $a>0$. Further we denote by $\zeta$ a random variable with the Pareto distribution defined above.


Finally, our random graph has the following parameters: 1) number of vertices~$n$, 2) vertex outdegree~$m$, 3) lifespan length~$N$, 4) quality exponent~$\gamma$, and 5) minimal quality~$a$.

\subsection{Degree distribution}

\subsubsection{Results}\label{sec:results}

In order to simplify the statements of theorems, we introduce an additional constant $\alpha$.
If $\gamma>2$, then we fix $\alpha = 2$.
If $1 < \gamma \le 2$, then $\alpha$ can be any constant such that $1< \alpha < \gamma$.

Let $N_n(d)$ be the number of vertices with degree $d$ in $G_n$.
We prove the following theorem.

\begin{thm}\label{th:expectation}
Assume that $d=d(n)$ increases with $n$ and $d=o\left( \min\left\{ \left(\frac{n}{N} \right)^{\frac{1}{\gamma+1}}\right.\right.$,  $\left.\left.N^{\frac{\alpha-1}{\gamma+\alpha+1}} \right\} \right)$, then
$$
\frac{\E [N_n(d)]}{n} = \frac{ \gamma }{d^{\gamma + 1}} \left( \frac{(\gamma-1)m}{\gamma} \right)^\gamma
                    \left( 1 + o(1) \right).
$$
\end{thm}

Theorem~\ref{th:expectation} shows that the expected number of vertices of degree $d$ decreases as $d^{-\gamma-1}$. In order to get the power-law degree distribution we also need to prove the concentration of the number of vertices with degree $d$ near its expectation.

\begin{thm}\label{th:concentration}
For every $d$ the following inequality holds:
$$
\Prob\left(|N_n(d) - \E [N_n(d)]| \ge  \sqrt{N n \log n}  \right) \le
 \frac{2}{\log n} .
$$
\end{thm}

Note that for $d = o\left(\left(\frac{n}{N\log{n}}\right)^{1/2(\gamma+1)}\right)$ we have $\sqrt{N n \log n} = o\left(n/d^{\gamma+1}\right)$,
so Theorem~\ref{th:concentration} gives the concentration.

We prove Theorem~\ref{th:expectation} in Sections~\ref{sec:weight} and~\ref{sec:expectation}. Theorem~\ref{th:concentration} is proven in Section~\ref{sec:concentration}.

\subsubsection{Concentration of the weight}\label{sec:weight}

Let us now fix some $n$ and $N=N(n)$. In this section we consider only the vertices $N \le p \le n-N+1$.

Let us denote by $\bar d(p)$ the degree of a vertex $p$ after its death and by $\bar d_{in}(p)$ the incoming degree of a vertex $p$ after its death, i.e., $\bar d_{in}(p) = \bar d(p) - m$.
By $Q(t)$ we denote the sum of qualities of the alive vertices at the $t$-th step, i.e.,
$$
Q(t) = \sum_{k=t-N}^{t-1} q(k)\,.
$$
We also say that $Q(t)$ is the weight of vertices at $t$-th step.
Note that

$$
\E \left[ \bar d_{in}(p) \mid q(p-N+1), \ldots, q(p+N-1) \right] = \sum_{i=1}^{N} \frac{m q(p)}{Q(p+i)} \,.
$$
Indeed, for each $1 \le i \le N$ the probability of an edge $(p+i,p)$ is equal to $\frac{mq(p)}{Q(p+i)}$ according to the definition of the model, since $Q(p+i)$ is the overall attractiveness of all vertices at $(p+i)$-th step.

Consider the lifespan of a vertex $p$ with a quality $q(p)$. We have $\E[Q(p+i)|q(p)] = q(p)+(N-1)\E[\zeta]$ for $1 \le i \le N$. We want to estimate the probability of this weight $\E[Q(p+i)|q(p)]$ to deviate from the value $N\E[\zeta]$.

Let $\xi_1, \dots \xi_{N-1}$ be the weights of vertices $p-N+1, \dots, p-1$ and $\eta_1, \dots \eta_{N-1}$ be the weights of vertices $p+1, \dots, p+N-1$. Let $W_p^q(i)$ be the overall weight of all living vertices when the age of $p$ equals $i$ given that $p$ has the quality $q$, i.e.,
$$
W_p^q(i) =  \sum_{k=1}^{i-1} \eta_k + q + \sum_{k=i}^{N-1} \xi_k \,.
$$

We will need the following auxiliary lemma.

\begin{lem}\label{lem1}
Let $\xi_1, \dots, \xi_n$ be mutually independent random variables, $\E [\xi_i] = 0$, $\E \left[|\xi_i|^{\alpha} \right] < \infty$, $1\le \alpha \le 2$, then
$$
\E \left[ |\xi_1 + \ldots + \xi_n|^{\alpha} \right] \le 2^{\alpha} \left(\E\left[|\xi_1|^{\alpha}\right]+\ldots+\E\left[|\xi_n|^{\alpha}\right] \right)\,.
$$
\end{lem}

We placed the proof of this lemma in the appendix.

\begin{thm}\label{Weight}
Consider a vertex $p$ with a quality $q$ such that $N \le p \le n-N+1$. If for some constant $c>0$ we have $|q-\E[\zeta]| \le N^c/2$, then
$$
\Prob\left(\max_{1\le i \le N}|W_p^q(i)-N\E[\zeta]| \ge N^c\right) = O\left( \frac{\E\left[|\zeta-\E[\zeta]|^{\alpha}\right]}{N^{\alpha c  -1}}\right)\,.
$$
\end{thm}

\begin{Proof}

Note that
\begin{multline*}
\Prob\left(\max_{1\le i \le N}|W_p^q(i)-\E[W_p^q(1)]| \ge x\right)
\\
\le \Prob\left(|W_p^q(1) -\E [W_p^q(1)] | \ge x/2\right) + \Prob\left(\max_{2\le i \le N} |W_p^q(i) -  W_p^q(1) | \ge x/2\right)\,.
\end{multline*}
Indeed, $\max_{1\le i \le N}|W_p^q(i)-\E [W_p^q(1)]| \le |W_p^q(1) -\E [W_p^q(1)] | + \max_{2\le i \le N} |W_p^q(i) -  W_p^q(1) |$ and if $\max_{1\le i \le N}|W_p^q(i)-\E [W_p^q(1)]| \ge x$ then either $|W_p^q(1) -\E [W_p^q(1)] | \ge x/2$ or $\max_{2\le i \le N} |W_p^q(i) -  W_p^q(1) | \ge x/2$.

In the case $\gamma > 2$ the random variables have finite variances and we can apply Chebyshev's and Kolmogorov's inequalities.

Chebyshev's inequality gives
$$
\Prob(|W_p^q(1)-\E [W_p^q(1)]| \ge x/2) \le \frac{4 N Var[\zeta]}{x^2}\,.
$$
Kolmogorov's inequality gives
\begin{multline*}
\Prob\left(\max_{2\le i \le N} |W_p^q(i) -  W_p^q(1) | \ge x/2\right)
\\
=
 \Prob\left(\max_{1\le i \le N-1}\left|\sum_{k=1}^i(\eta_{k}-\xi_{k})\right| \ge x/2\right) \le \frac{8 N \, Var[\zeta]}{x^2}\,.
\end{multline*}
So, finally we get
$$
\Prob\left(\max_{1\le i \le N}|W_p^q(i)- \E [W_p^q(1)]| \ge x\right) \le \frac{12 N \, Var[\zeta]}{x^2}\,.
$$
Take $x = N^{c}/2$ and note that $|\E[ W_p^q(1)] - N\E [\zeta]| = |q-\E[\zeta]| \le N^c/2$. Therefore we get
\begin{multline*}
\Prob\left(\max_{1\le i \le N}|W_p^q(i)-N\E[\zeta]| \ge N^c\right)
\le
 \Prob\left(\max_{1\le i \le N}|W_p^q(i)- \E [W_p^q(1)]| \ge N^c/2\right)
 \\
 \le
 \frac{48 \, Var[\zeta]}{N^{2c-1}}
=
O\left( \frac{\E\left[|\zeta-\E[\zeta]|^{\alpha}\right]}{N^{\alpha c  -1}}\right)\,,
\end{multline*}
since $\alpha = 2$ in this case.

Now consider the case $1<\gamma\le2$.
In this case we have $1 < \alpha < \gamma$.
We cannot use Chebyshev's inequality now, but we can apply Markov's inequality and Lemma~\ref{lem1}:
\begin{multline*}
\Prob(|W_p^q(1)-\E [W_p^q(1)]| \ge x/2) =
\Prob\left(|W_p^q(1)-\E [W_p^q(1)]|^{\alpha}\ge (x/2)^{\alpha}\right)
\\
\le\frac{\E\left[|W_p^q(1)-\E[W_p^q(1)]|^{\alpha}\right]}{(x/2)^{\alpha}} \le
 \frac{4^{\alpha}N\E\left[|\zeta-\E[\zeta]|^{\alpha}\right]}{x^{\alpha}}\,.
\end{multline*}

Instead of Kolmogorov's inequality, we use Doob's martingale inequality and Lemma~\ref{lem1}.
Note that $S_i = \left|\sum_{j=1}^i(\eta_{j}-\xi_{j})\right|$ is a submartingale as a convex function of a martingale. Thus,
\begin{multline*}
\Prob\left(\max_{1\le i \le N-1}\left|\sum_{j=1}^i(\eta_{j}-\xi_{j})\right| \ge x/2\right) \le \frac{\E\left[\left|\sum_{j=1}^{N-1}(\eta_{j}-\xi_{j})\right|^{\alpha}\right]}{(x/2)^{\alpha}}
\\
\le \frac{4^{\alpha}N\E\left[\left|\eta_{1}-\xi_{1}\right|^{\alpha}\right]}{x^{\alpha}} \le \frac{8^{\alpha}N\E\left[|\zeta-\E[\zeta]|^{\alpha}\right]}{x^{\alpha}}\,.
\end{multline*}
So, finally we get
$$
\Prob\left(\max_{1\le i \le N}|W_p^q(i)- \E[W_p^q(1)]| \ge x\right) \le \frac{4^{\alpha}\left(2^{\alpha}+1\right) N \, \E\left[|\zeta-\E[\zeta]|^{\alpha}\right]}{x^{\alpha}}\,.
$$
Now take $x = N^{c}/2$ and note that $|\E [W_p^q(1)] - N\E [\zeta]| \le N^c/2$. As before, we can estimate
\begin{multline*}
\Prob\left(\max_{1\le i \le N}|W_p^q(i)-N\E[\zeta]| \ge N^c\right) \le
\frac{8^{\alpha}\left(2^{\alpha}+1\right) N \, \E\left[|\zeta-\E[\zeta]|^{\alpha}\right]}{N^{\alpha c}} \\ \le
\frac{320 \,\E\left[|\zeta-\E[\zeta]|^{\alpha}\right]}{N^{\alpha c -1}}
=
O\left( \frac{\E\left[|\zeta-\E[\zeta]|^{\alpha}\right]}{N^{\alpha c  -1}}\right)\,.
\end{multline*}

This concludes the proof of Theorem~\ref{Weight}.

\end{Proof}


\subsubsection{Expectation}\label{sec:expectation}

Let $\rho(d,q)$ be the conditional probability that a vertex $p$ such that $N \le p \le n-N+1$ with a quality $q$ has an in-degree $d$, i.e., $\rho(d,q)= \Prob(\bar d_{in}(p) = d | q(p)=q)$.
Note that $\rho(d,q)$ does not depend on $p$.
By $N_n^{in}(d)$ we denote the number of vertices with in-degree $d$, so $N_n^{in}(d) = N_n(d+m)$.
The expectation of $N_n^{in}(d)$ is
\begin{equation}\label{eq:expectation}
\E \left[ N_n^{in}(d) \right] = (n-2N) \int_{a}^{\infty} f(q) \rho(d,q) d q + r(N)\,,
\end{equation}
where $f(q)$ is the density function of Pareto distribution defined above and $r(N)$, $0 \le r(N) \le 2N$, is the error term. We have this error term since the first and the last $N$ vertices behave differently.

Let $c$ be some positive constant. We estimate the integral
$$
I = \int_a^\infty f(q) \rho(d,q) d q =
\int_{a}^{N^c/2} f(q) \rho(d,q) d q + \int_{N^c/2}^\infty f(q) \rho(d,q) d q = I_1 + I_2.
$$

Note that
\begin{equation}\label{eq:I_2}
I_2 = \int_{N^c/2}^\infty f(q) \rho(d,q) d q \le \int_{N^c/2}^\infty f(q) d q = \int_{N^c/2}^\infty \frac{\gamma  a^{\gamma}}{q^{\gamma+1}} d q =  \frac{(2a)^{\gamma}}{N^{c\gamma}}\,.
\end{equation}

Consider the event
$$A = { \left\{\max_{1\le i \le N}|Q(p + i)- N \E[\zeta]| \le N^c\right\} }$$
and the following conditional probabilities:
$$\rho_A(d,q) = \Prob(\bar d_{in}(p) = d | q(p)=q, A)\,,$$
$$\rho_{\bar A}(d,q) = \Prob(\bar d_{in}(p) = d | q(p)=q, \bar A)\,.$$
Then we have
\begin{equation}\label{eq:rho}
 \rho(d,q) = \rho_A(d,q) \Prob(A | q(p)=q) + \rho_{\bar A}(d,q) \Prob(\bar A | q(p)=q)\,.
\end{equation}

Let us use this representation to split $I_1$ into two integrals using \eqref{eq:rho}:
$$
I_1 = \int_{a}^{N^c/2} f(q) \rho_A(d,q) \Prob(A_q) d q + \int_{a}^{N^c/2} f(q)  \rho_{\bar A}(d,q) \Prob(\bar A_q) d q = I_1^1 + I_1^2\,,
$$
where we use the following notation:

$$A_q = [ A | q(p)=q ] = { \left\{\max_{1\le i \le N}| W_p^q(i) - N \E[\zeta]| \le N^c\right\} } \,, $$
$$\bar A_q = [ \bar A | q(p)=q ] = { \left\{\max_{1\le i \le N}| W_p^q(i) - N \E[\zeta]| > N^c\right\} } \, .$$

Let us assume that $N^c/2 > \E[\zeta]$, this holds if $N$ is large enough (the fact that $N$ grows follows from the statement of Theorem~\ref{th:expectation}, while $\E[\zeta]$ is constant). Note that
\begin{equation}\label{eq:I_1^2}
I_1^2 \leq \max_{q \le N^c/2} \Prob(\bar A_q)
\end{equation}
and since $q \le N^c/2$ Theorem~\ref{Weight} gives us an upper bound for it, i.e.,
\begin{equation}\label{eq:2}
\max_{q \le N^c/2} \Prob(\bar A_q) = O\left( N^{1 - \alpha c} \right)\,, 
\end{equation}
where $\alpha = 2$ for $\gamma > 2$ and
 $1 < \alpha < \gamma$ for $1<\gamma\le2$.

So, let us now focus on $I_1^1$. First we estimate $\rho_A(d,q)$. Recall that $\rho_A(d,q) = \Prob(\bar d_{in}(p) = d | q(p)=q, A)$.
Note that during the lifespan of a vertex $p$ there are $mN$ mutually independent edges which may lead to $p$.
For an edge from a vertex $p+i$ the probability to choose $p$ is $\frac{q}{W_p^q(i)}$.
Given the event $A_q$ we have $N\E[\zeta] - N^c \le W_p^q(i) \le N\E[\zeta] + N^c$.
Therefore we have the following bounds for $\rho_A(q, d)$:
\begin{multline*}
{mN \choose d}\left(\frac{q}{N\E[\zeta]+N^c}\right)^{d} \left(1 - \frac{q}{N\E[\zeta]-N^c}\right)^{mN-d}
\\
\le \rho_A(d,q) \le
\\
\le {mN \choose d}\left(\frac{q}{N\E[\zeta]-N^c}\right)^{d} \left(1 - \frac{d}{N\E
[\zeta]+N^c}\right)^{mN-d}.
\end{multline*}

Thus, $\left( 1 - \max_{q \le N^c/2} \Prob(\bar A_q) \right) S_- \leq I_1^1 \leq S_+$ where

$$
S_{\mp} = \int_{a}^{N^c/2} f(q) {mN \choose d} \left(\frac{q}{N\E[\zeta] \pm N^c}\right)^{d} \left(1 - \frac{q}{N\E[\zeta] \mp N^c}\right)^{mN-d} d q\,.
$$

We will use the following lemma.

\begin{lem}\label{S-+}
Assume that both $N$ and $d$ grow but ${d = o(N^{1-c})}$. If $1/2 \le c \le 1$, then
$$
S_{\mp} =  \frac{ \gamma }{d^{\gamma + 1}} \left( \frac{(\gamma-1)m}{\gamma} \right)^\gamma
           \left( 1 + o(1) \right)\,.
$$
\end{lem}

We placed the proof of this technical lemma in the appendix.
Now we will use this lemma to prove the theorem. Using Equations \eqref{eq:expectation}, \eqref{eq:I_2}, and \eqref{eq:I_1^2} we get the following bounds for $\frac{\E [N_n(d+m)]}{n}$:
\begin{multline}\label{eq:bounds}
\left(1-\frac{2N}{n}\right)\underbrace{\left( 1 - \max_{q \le N^c/2} \Prob(\bar A_q) \right) S_- }_{\leq I_1^1}
\leq  \frac{\E [N_n(d+m)]}{n}
\\
\le\underbrace{S_+}_{\geq I^1_1} +
\underbrace{\max_{q \le N^c/2} \Prob(\bar A_q)}_{\geq I^2_1} +
\underbrace{\frac{(2a)^{\gamma}}{N^{c\gamma}}}_{\geq I_2} + \frac{2N}{n}\,.
\end{multline}

Now we show that for some parameter $c$ all error terms in Equation \eqref{eq:bounds} are negligible in comparison with the main term $d^{-\gamma-1}$ from Lemma~\ref{S-+}.
We take $c = \frac{\gamma+2}{\gamma+\alpha+1}$. Note that we can apply Lemma~\ref{S-+} since ${d = o(N^{1-c})}$ due to the statement of Theorem~\ref{th:expectation}.

\begin{itemize}
\item[1.] $\frac{N}{n} = o\left(d^{-\gamma-1}\right)$, since $d=o\left( \left(\frac{n}{N} \right)^{\frac{1}{\gamma+1}} \right)$.
\item[2.] $\frac{(2a)^{\gamma}}{N^{c\gamma}} = o(d^{-\gamma-1})$ if $d = o(N^{c\gamma/(\gamma+1)}).$ This holds for $c = \frac{\gamma+2}{\gamma+\alpha+1}$ and $d = o \left( N^{\frac{\alpha-1}{\gamma+\alpha+1}} \right)$.
\item[3.] For $\gamma \le 2$, $\Prob\left(\bar A_q\right) = O\left( N^{1-c\alpha} \right) = o(d^{-\gamma-1})$ if $d = o(N^{(c\alpha-1)/(\gamma+1)})$, i.e., $d = o\left(N^{\frac{\alpha-1}{\gamma+\alpha+1}}\right)$. Here we used Equation~\eqref{eq:2}.
\end{itemize}

We demonstrated that all the error terms in Equation~\eqref{eq:bounds} equal $o(d^{-\gamma-1})$.
Therefore, from Lemma~\ref{S-+} we get
$$
\frac{\E [N_n(d+m)]}{n} = \frac{ \gamma }{d^{\gamma + 1}} \left( \frac{(\gamma-1)m}{\gamma} \right)^\gamma
           \left( 1 + o(1) \right)\,.
$$

To conclude the proof of Theorem~\ref{th:expectation} it remains to note that the asymptotic for $N_n(d)$ is the same as for $N_n(d+m)$.

\subsubsection{Concentration}\label{sec:concentration}

We use Chebyshev's inequality to prove concentration. In order to do this we first estimate $Var[N_n(d)]$. Note that if $|i-j| \ge N$ then the degrees of $i$ and $j$ are independent. Therefore
\begin{multline*}
Var[N_n(d)] \\= \sum_{i,j=1}^n \left( \Prob (d_n(i) = d, d_n(j) = d) - \Prob (d_n(i) = d)\Prob(d_n(j) = d)\right) \le 2nN\,.
\end{multline*}

Applying Chebyshev's inequality we get
$$
\Prob\left(|N_n(d) - \E [N_n(d)]| > \sqrt{Nn\log n}\right) \le \frac{Var[N_n(d)]}{Nn\log n} \le \frac{2}{\log n}\,.
$$

\textbf{Remark.}
\textit{Note that instead we could use Azuma--Hoeffding inequality, since  $|\E[N_n(d)|G_i] - \E[N_n(d)|G_{i-1}]| \le (N + 1) m$. In this case we get
$$
\Prob\left(|N_n(d) - \E [N_n(d)]| \ge  \sqrt{n \log n} (N+1) \right) \le 2  n^{-1/2 m^2} .
$$
So, on the one hand the range of degrees for which we get concentration is smaller in this case. We get concentration for $d = o\left(\left(\frac{\sqrt{n}}{N\sqrt{\log{n}}}\right)^{1/(\gamma+1)}\right)$. On the other hand, the concentration is tighter, so we can say that for all $d$ in this range the number of vertices of degree $d$ is near its expectation.}

\subsection{Recency property}\label{sec:recency}

Let $e(T)$ be the fraction of edges in a graph which connect vertices with age difference greater than $T$, i.e., vertices $i$ and $j$ with $|i-j|>T$.
In \cite{2012model} a \textit{recency property} was introduced, which reflects the fact that new media pages tend to connect to other media pages of similar age.
Namely, for the media related part of the Web it was shown that $e(T)$ decreases exponentially fast. In this section we show that we have linear decay of $e(T)$ for the model under consideration.

\begin{thm}
For any integer $T$
$$\E [e(T)] =
\begin{cases}
1 - \frac{T}{N} + O\left(\frac{N}{n}\right),&\text{if $T\le N$;}\\
0,&\text{if $T > N$.}
\end{cases}
$$
\end{thm}

\begin{Proof}
Consider any vertex $n>N$ and any edge $ni$, $i<n$, drawn from this vertex. The probability that $n-i > T$ is the probability to choose one vertex from $n-N, \dots, n-T-1$. Since qualities of vertices are i.i.d. random variables, this probability equals $\frac{N-T}{N}$. From this the theorem follows.
\end{Proof}

\begin{thm}
$$
\Prob\left(|e(T) - \E [e(T)]| \ge  \sqrt\frac{N\log n}{n} \right) \le \frac{2}{\log n} .
$$
\end{thm}

\begin{Proof}
Here we again use Chebyshev's inequality. Let $e_i$ and $e_j$ be any two different edges drawn from vertices $i$ and $j$. For an edge $e$ let $l(e)$ be the age difference between the endpoints of the edge.

Note that if $|i-j| \ge N$, then $l(e_i)$ and $l(e_j)$ are independent:
$$
\Prob(l(e_i) > T, l(e_j) > T) - \Prob(l(e_i) > T)\Prob(l(e_j) > T) = 0\,.
$$



From this we get $Var[m\,n\,e(T)] \le 2 m^2 N n$, since we take into account only the summands corresponding to edges $e_i$ and $e_j$ with $|i-j| \le N$ and $\Prob(l(e_i) > T, l(e_j) > T) - \Prob(l(e_i) > T)\Prob(l(e_j) > T) \le 1$.

Therefore,
$$
\Prob\left(m\,n\,|e(T) - \E [e(T)]| \ge  m \sqrt{N \,n \log n} \right) \le \frac{Var[m\,n\,e(T)]}{m^2\, N\,n \log n} \le \frac{2}{\log n}\,.
$$

\end{Proof}


Finally, let us discuss an intuition behind the recency factor $I[i > t - N]$.
This factor has the following natural interpretation.
Links to a lot of media pages can usually be found on some pages which are content sources.
And new pages are popular while they can be found on such content sources.
After some period of time other new pages appear on a content source and they replace old ones. Therefore, it seems natural to assume that after some period of time old pages become unpopular.
However, as it turned out, this recency factor can model only linear decay of $e(T)$, while we want to see exponential one.
One possible reason for this lack of agreement is that the most recent links are usually placed near the top of a page and they may attract more attention because of this.

\section{Attractiveness function $q(i) e^{-\frac {t-i}{N}}$}

Now we switch to the attractiveness function $q(i) e^{-\frac {t-i}{N}}$.
In this case, the popularity of a vertex decreases exponentially with the age of the vertex.
Again, we assume that the random variables $\zeta_1, \zeta_2, \ldots$ have the Pareto distribution with the density function $f(x) = \frac{\gamma  a^{\gamma} I[x > a]}{x^{\gamma+1}}$, where $\gamma > 1$, $a>0$. And $\zeta$ again is a random variable with the Pareto distribution defined above.

\subsection{Degree distribution}

\subsubsection{Results}

For the model with the exponential recency factor we get the results similar to ones for the model with the indicator recency factor (see Section~\ref{sec:results}).

Recall that a constant $\alpha$ is defined as follows: if $\gamma>2$, then $\alpha = 2$;
if $1 < \gamma \le 2$, then $\alpha$ can be any constant such that $1< \alpha < \gamma$.

\begin{thm}\label{th:expectation2}
If $d=d(n)$ increases with $n$ and
$d=o\left( \min\left\{ \left(\frac{n}{N\log N} \right)^{\frac{1}{\gamma+1}} \right.\right.$,
$\left.\left. N^{\frac{\alpha-1}{\alpha+(\gamma+1)(\alpha+1)}} \right\}\right)$, then
$$
\frac{\E [N_n(d)]}{n} = \frac{ \gamma }{d^{\gamma + 1}} \left( \frac{(\gamma-1)m}{\gamma} \right)^\gamma
                    \left( 1 + o(1) \right)\,.
$$
\end{thm}
Again, the expectation of the number of vertices with degree $d$ decreases as $d^{-\gamma-1}$. The next theorem shows that the number of vertices of degree $d$ is concentrated near its expectation.

\begin{thm}\label{th:concentration2}
For every $d$ the following inequality holds:
$$
\Prob\left(|N_n(d) - \E [N_n(d)]| > \sqrt{Nn\log n}\right) = O\left(\frac{1}{\log n}\right)\,.
$$
\end{thm}
As before, for $d = o\left(\left(\frac{n}{N\log{n}}\right)^{\frac{1}{2(\gamma+1)}}\right)$ we have $\sqrt{N n \log n} = o\left(\E [N_n(d)]\right)$ and Theorem~\ref{th:concentration2} gives the concentration.

We prove Theorem~\ref{th:expectation2} in Sections~\ref{sec:weight2} and~\ref{sec:expectation2}. Theorem~\ref{th:concentration2} is proven in Section~\ref{sec:concentration2}.

\subsubsection{Concentration of the overall attractiveness}\label{sec:weight2}

We fix some $n$ and $N=N(n)$.

By $Q(t)$ we denote the total attractiveness of all vertices at $t$-th step, i.e.,
$$
Q(t) = \sum_{k=1}^{t-1} q(k)e^{-\frac{t-k-1}{N}}\,.
$$
The average value of $Q(t)$ is
$$
\E [Q(t)] = \E [\zeta] \sum_{k=0}^{t-2} e^{-\frac{k}{N}} = \E [\zeta] \frac{1 - e^{-\frac{t-1}{N}}}{1 - e^{-\frac{1}{N}}}
= N \E [\zeta] \left(1 + O\left( e^{-t/N} \right) + O\left(1/N\right)  \right)\,.
$$
If $t > N \log N$, then 
$$
\E [Q(t)] = N \E [\zeta] \, (1 + O(1/N))\,.
$$

Again, by $W_p^q(i)$ we denote the total attractiveness of all vertices when the age of $p$ equals $i$ given the quality $q$ of the vertex $p$.

\begin{thm}\label{Weight2}
Fix some positive constant $c$. Let $\varphi(N)$ be any function such that $\varphi(N) > \log(CN)$ for some $C>0$. Then for any $p>N \varphi(N)$
with a quality $q$ satisfying $ { \mid q - \E [\zeta] \mid \le N^c/3} $ we have:\\
$$
\Prob\left(\max_{1\le i \le N\varphi(N)}|W_p^q(i)-N\E [\zeta]| \ge N^c\right) = O\left( e^{\alpha\varphi(N)}N^{1-\alpha c} \right).
$$
\end{thm}

\begin{Proof}


Note that $\E \left[ Q(p + i + 1) \mid Q(p+i) \right] = Q(p + i) e^{-\frac{1}{N}} + \E [\zeta] $.
Therefore $X_i = e^{\frac{i}{N}} \cdot \left( Q(p + i)  - \frac{\E [\zeta]}{1 - e^{-\frac{1}{N}}} \right)$ is a martingale. Indeed,
$$
\E \left[ X_{i+1} \mid X_i \right]
= \left( Q(p + i) e^{-\frac{1}{N}} + \E [\zeta] \right) e^{\frac{i+1}{N}} - \frac{e^{\frac{i+1}{N}}\E [\zeta]}{1 - e^{-\frac{1}{N}}}
= X_i \,.
$$
So, we can apply Doob's inequality for a submartingale $|X_i|$:
\begin{multline*}
\Prob\left(\max_{1\le i \le N \varphi(N)}\left|e^{\frac{i}{N}}\left(Q(p + i) - \frac{\E [\zeta]}{1-e^{-1/N}}\right)\right| \ge x\right)
\\
\le \frac{\E\left[\left|e^{\frac{N\varphi(N)}{N}}\left(Q(p + N\varphi(N)) - \frac{\E [\zeta]}{1-e^{-1/N}}\right)\right|^{\alpha}\right]}{x^{\alpha}}\,.
\end{multline*}
So, we get
\begin{multline*}
\Prob\left(\max_{1\le i \le N\varphi(N)}\left|Q(p+i) - \frac{\E [\zeta]}{1-e^{-1/N}}\right| \ge N^{c}/3\right)
\\
\le \frac{3^\alpha e^{\alpha\varphi(N)}\E\left[\left|Q(p + N\varphi(N)) - \frac{\E [\zeta]}{1-e^{-1/N}}\right|^{\alpha}\right]}{N^{\alpha c}}\,.
\end{multline*}

First, consider the case $\gamma > 2$. In this case we have $\alpha = 2$
Using
$$
\frac{\E [\zeta]}{1-e^{-1/N}} = \E [Q(p + N\varphi(N))] + \E [\zeta] \sum_{k=p + N\varphi(N)-1}^{\infty} e^{-\frac k N}
$$
we get
\begin{multline*}
\E \left[\left(Q(p + N\varphi(N)) - \frac{\E \zeta}{1-e^{-1/N}}\right)^{2}\right]
\\
= \E\left[\left(Q(p + N\varphi(N)) - \E Q(p + N\varphi(N))\right)^2 \right] +  (\E [\zeta])^2 \left( \frac{e^{-\frac{p + N\varphi(N)-1}{N}}}{1-e^{-1/N}} \right)^2
\\
= Var[\zeta]\sum_{k=0}^{p + N\varphi(N)-2} e^{-\frac{2k}{N}}  + O \left( \frac{e^{-\frac{2(N \log (C N) + N\varphi(N))}{N}}}{(1-e^{-1/N})^2} \right)
\\
= O\left(\frac{1}{1-e^{-2/N}} + e^{-2\varphi(N)}\right)
= O \left(N \right)\,.
\end{multline*}
So,
$$
\Prob\left(\max_{1\le i \le N\varphi(N)}\left|Q(p+i) - \frac{\E [\zeta]}{1-e^{-1/N}}\right| \ge N^{c}/3\right) = O\left( e^{2\varphi(N)}N^{1-2c} \right)\,.
$$
Now we can estimate $W_p^q(i)$ which is $Q(p+i)$ given the quality $q$ of the vertex $p$.
We have ${\mid q - \E [\zeta] \mid \le N^c/3}$ and $\left| \frac{\E [\zeta]}{1-e^{-1/N}} - N\E [\zeta]\right| \le N^c/3$ for large $N$, therefore
$$
\Prob\left(\max_{1\le i \le N\varphi(N)}|W_p^q(i)- N\E [\zeta]| \ge N^c\right)
= O\left( e^{2\varphi(N)}N^{1-2c} \right)\,.
$$

Similarly, for $\gamma \le 2$ using Lemma~\ref{lem1} we get
$$
\Prob\left(\max_{1\le i \le N\varphi(N)}|W_p^q(i)- N\E [\zeta]| \ge N^c\right) = O\left( e^{\alpha\varphi(N)}N^{1-\alpha c} \right)\,.
$$

\end{Proof}

\subsubsection{Expectation}\label{sec:expectation2}

Let $\varphi(N)$ be any function such that $\varphi(N) > \log(CN)$ for some $C>0$.
Let $\rho(d,q)$ be the conditional probability that a vertex $p$ such that  $N\varphi(N) \le p \le n - N \varphi(N)+1$ has an in-degree $d$ given a quality $q$ of this vertex, i.e., $\rho(d,q)= \Prob(\bar d_{in}(p) = d | q(p)=q)$.  We omit $n$ and $p$ in the notation $\rho(d,q)$ because, as we will see, we get similar bounds for $\rho(d,q)$ for all $p$ such that $N\varphi(N) \le p \le n - N \varphi(N)+1$. Using this notation, we get the following equality:

\begin{equation}\label{eq:1'}
\E[ N_n^{in}(d) ] = (n-2N\varphi(N)) \int_{a}^{\infty} f(q) \rho(d,q) d q + r(N)\,,
\end{equation}
where $f(q)$ is the density function of Pareto distribution and $r(N)$, $0 \le r(N) \le 2N\varphi(N)$ is the error term.

Let $r$ and $c$ be some constants such that $0<r<1/2$ and $1/2<c<1$.
As in Section~\ref{sec:expectation}, we split the integral
$$
I = \int_{a}^{\infty} f(q) \rho(d,q) d q =
\int_{a}^{N^r} f(q) \rho(d,q) d q + \int_{N^r}^{\infty} f(q) \rho(d,q) d q = I_1 + I_2
$$
and
\begin{equation}\label{eq:I_2'}
I_2 \le \int_{N^r}^{\infty} f(q) d q =  \frac{a^{\gamma}}{N^{r\gamma}}\,.
\end{equation}

The event $A$ is defined as in Section~\ref{sec:expectation}:
$$A = { \left\{\max_{1\le i \le N \varphi(N)}|Q(p + i)- N \E[\zeta]| \le N^c\right\} }\,.$$

We again split $I_1$ into two integrals:
$$
I_1 = \int_{a}^{N^r} f(q) \rho_A(d,q) \Prob(A_q) d q + \int_{a}^{N^r} f(q)  \rho_{\bar A}(d,q) \Prob(\bar A_q) d q = I_1^1 + I_1^2 \,,
$$
where
$$A_q = [ A | q(p)=q ] = { \left\{\max_{1\le i \le N\varphi(N)}| W_p^q(i) - N \E[\zeta]| \le N^c\right\} } \,,$$
$$\bar A_q = [ \bar A | q(p)=q ] = { \left\{\max_{1\le i \le N\varphi(N)}| W_p^q(i) - N \E[\zeta]| > N^c\right\} } \,.$$

We can estimate
\begin{equation}\label{eq:I_1^2'}
I_1^2 \leq \max_{q \leq N^r} \Prob(\bar A_q)
\end{equation}
 and for $q \le N^r$
Theorem~\ref{Weight2} gives the upper bound for $\Prob(\bar A_q)$ (since $N^r < N^c/3$ and $|q - \E [\zeta]| \le N^r$ for large $N$):
\begin{equation}\label{eq:max1}
\max_{q \leq N^r} \Prob(\bar A_q) = O\left( e^{\alpha\varphi(N)}N^{1-\alpha c} \right)  \,.
\end{equation}

Let us now focus on $I_1^1$.
Consider an event $R_p^q(k)$ that there is an edge from at least one vertex $p+i$ with $i \geq k$ to a vertex $p$ with a quality $q$.
Then for $k>N$ conditional probability of $R_p^q(k)$ given $A_q$ can be estimated as follows
\begin{multline*}
\Prob\left(R_p^q(k) \mid A_q \right) \le \sum_{i=k}^{\infty} \Prob\left( \text{edge } (p+i, p) \text{ belongs to } G_n \mid A_q \right)
\\
 \le  \sum_{i=k}^{\infty} \frac{m\, q\, e^{-\frac{i-1}{N}}}{ \sum_{j=0}^{N} a e^{\frac{-j}{N}}}
 \le  \sum_{i=k}^{\infty} \frac{m\, q\, e^{-\frac{i-1}{N}}}{aN /2} = O\left(q e^{-\frac{k}{N}}\right)\,.
\end{multline*}
This estimate means that the most contribution to the final degree of a vertex is made during
the first several steps after its appearance and we have the following bounds for $\rho_A(d,q)$:
$$
\rho_A(d,q) = \rho_{\mp}(d,q) + O \left(q e^{-\varphi(N)}\right),
$$
where $\rho_{\mp}(d,q)$ are lower and upper bounds for the probability that a vertex $p$ with a quality $q$ has an in-degree $d$ in $\tilde{G}_{p+N\varphi(N)}^n$ given $A_q$. We can estimate $\rho_{\mp}(d,q)$ in the following way. A vertex $p$ has an in-degree $d$ in $\tilde{G}_{p+N\varphi(N)}^n$ if $d$ edges out of $m \varphi(N) N$ are connected to this vertex and others are not. For every set of indexes $0 \leq i_1 < \ldots < i_d \leq m \varphi(N) N$ we should multiply the probabilities that the corresponding edges go to the vertex $p$. Given $A_q$, these probabilities can be estimated by $\frac{q e^{\frac{-[ i_j / m ]}{N}}}{N \E [\zeta] \mp N^c}$. And we should also multiply the obtained product by the probabilities that other edges are not connected to $p$, i.e., $\left( 1 - \frac{q e^{\frac{-[ i / m ]}{N}}}{N \E [\zeta] \pm N^c}  \right)$ for the corresponding indexes $i$. Finally, we get:

$$
\rho_{\mp}(d,q) =
\prod_{i=0}^{m \varphi(N) N} \left( 1 - \frac{q e^{\frac{-[ i / m ]}{N}}}{N \E [\zeta] \pm N^c}  \right)
\sum_{0 \leq i_1 < \ldots < i_d \leq m \varphi(N) N}
\prod_{j=1}^{d} \frac{ \frac{q e^{\frac{-[ i_j / m ]}{N}}}{N \E [\zeta] \mp N^c}}
                     {1 - \frac{q e^{\frac{-[ i_j / m ]}{N}}}{N \E [\zeta] \pm N^c}} \,.
$$

Now we put
$$
S_{\mp}(d,q) := \int_{a}^{N^r} f(q) \rho_{\mp}(d,q) d q.
$$
Using this notation, we can estimate $I_1^1$ in the following way:
\begin{equation}\label{eq:I_1^1u}
I_1^1 \le \int_{a}^{N^r} f(q) \rho_A(d,q) d q \le S_{+} + O\left( \int_a^\infty f(q)  q e^{-\varphi(N)} dq  \right)\,,
\end{equation}
\begin{equation}\label{eq:I_1^1l}
I_1^1 \ge \left( 1 - \max_{q\le N^r} \Prob(\bar A_q) \right) S_- + O\left( \int_a^\infty f(q)  q e^{-\varphi(N)} dq  \right) \,.
\end{equation}

We estimate $S_{\mp}$ in the following way.

\begin{lem}\label{rho-+_exp}

Assume that both $d$ and $N$ grow, $d = o(N^{1-c})$, $d = o(e^{\varphi(N)})$, and $q \leq N^r$, then
$$
S_{\mp}(d,q) = \frac{ \gamma }{d^{\gamma + 1}} \left( \frac{(\gamma-1)m}{\gamma} \right)^\gamma \left(1 + o(1)\right)\,.
$$
\end{lem}

We placed the proof of this technical lemma in the appendix.

Finally, using Equations \eqref{eq:1'}, \eqref{eq:I_2'}, \eqref{eq:I_1^2'}, \eqref{eq:I_1^1u}, and \eqref{eq:I_1^1l}, we get
\begin{multline}\label{eq:final}
\left(1-\frac{2N\varphi(N)}{n}\right)\underbrace{\left( 1 - \max_{q\le N^r} \Prob(\bar A_q) \right) S_-+ O\left( \int_a^\infty f(q)  q e^{-\varphi(N)} dq  \right)}_{\leq I_1^1}
\\
\le \frac{\E [N_n(d+m)]}{n}
\\
\leq
\underbrace{S_+ + O\left( \int_a^\infty f(q)  q e^{-\varphi(N)} dq  \right)}_{\geq I_1^1}
+ \underbrace{\max_{q\le N^r} \Prob(\bar A_q)}_{\geq I_1^2}
+ \underbrace{\frac{a^{\gamma}}{N^{r\gamma}}}_{\geq I_2} + \frac{2N\varphi(N)}{n}\,.
\end{multline}

We want all the error terms in Equation~\eqref{eq:final} to be $o \left( d^{-\gamma-1} \right)$. In order to do this, we need
to find the proper values of $c$ and $\varphi(n)$.
Note that we have already assumed that ${d = o\left(N^{1-c}\right)}$ and $d = o\left(e^{\varphi(N)}\right)$.
We have to show that the following conditions hold.
\begin{enumerate}
\item $N\varphi(N)/n = o(d^{-\gamma-1})$ if $d = o\left(\left(\frac{n}{N\varphi(N)}\right)^{\frac{1}{\gamma+1}}\right)$\,. This holds under the conditions of the theorem.
\item $\frac{a^{\gamma}}{N^{r\gamma}} = o(d^{-\gamma-1})$ if $d = o(N^{r\gamma/(\gamma+1)})$. Put $r = \frac{9}{22}$, then we have $d = o(N^{r\gamma/(\gamma+1)})$ under the conditions of the theorem since $N^{9\gamma/22(\gamma+1)} \ge N^{\frac{\alpha-1}{\alpha+(\gamma+1)(\alpha+1)}}$ for both $\gamma > 2$ and $1 \le \gamma \le 2$.
\item
$\max_{q\le N^r} \Prob(\bar A_q) = O\left(e^{\alpha\varphi(N)}N^{1-\alpha c}\right) = o(d^{-\gamma-1})$ if \\ $e^{\alpha\varphi(N)} = o\left(N^{\alpha c-1} d^{-\gamma-1}\right)$. Here we used Equation~\eqref{eq:max1}.
\item $\int_a^\infty f(q) O\left(q e^{-\varphi(N)} \right) dq = O \left( e^{-\varphi(N)} \right) = o \left( d^{-\gamma-1} \right)$ if $d^{\gamma+1} = o \left( e^\varphi(N) \right)$\,.
\end{enumerate}

We take $c = \frac{1+(\gamma+1)(\alpha+1)}{\alpha+(\gamma+1)(\alpha+1)}$ and $\varphi(N) = \log{\frac{N(\alpha-1)(\gamma+1)}{\alpha+(\gamma+1)(\alpha+1)}}$. Then for
$d = o \left( N^{\frac{\alpha-1}{\alpha+(\gamma+1)(\alpha+1)}} \right)$ all the above conditions hold.
This means that all the error terms in Equation~\eqref{eq:final} equal $o(d^{-\gamma-1})$.
Therefore, we obtained the required asymptotic for  $\frac{\E [N_n(d+m)]}{n}$.

To conclude the proof of Theorem~\ref{th:expectation2} it remains to note that the asymptotic for $N_n(d)$ is the same as for $N_n(d+m)$.

\subsubsection{Concentration}\label{sec:concentration2}

We prove Theorem~\ref{th:concentration2} using Chebyshev's inequality. In order to apply this inequality we first estimate $Var[N_n(d)]$:
$$
Var[N_n(d)] = \sum_{i,j=1}^n \left( \Prob (d_n(i) = d, d_n(j) = d) - \Prob (d_n(i) = d)\Prob(d_n(j) = d) \right)\,.
$$

Let us estimate the difference $\Prob (d_n(i) = d, d_n(j) = d) - \Prob (d_n(i) = d)\Prob(d_n(j) = d)$ for $i < j$.

Note that
\begin{equation}\label{eq:independence}
\Prob (d_{j}(i) = d, d_n(j) = d) = \Prob (d_{j}(i) = d)\Prob(d_n(j) = d)\,.
\end{equation}
In order to prove this we first show that \eqref{eq:independence} holds given all the qualities $q_1, \ldots, q_n$ and then integrate over all qualities. Given the qualities, $\Prob (d_{j}(i) = d, d_n(j) = d)$ is the sum over all $mi<i_1<\ldots<i_d\le mj$, $mj<j_1<\ldots<j_d\le mn$ of the probabilities that the corresponding edges $([i_k/m], i)$ and $([j_k/m], j)$ are drawn and all other edges $(i', i)$ with $i<i'\le j$ and $(j',j)$ with $j<j'\le n$ are absent. Since qualities are fixed, these events are independent and $\Prob (d_{j}(i) = d, d_n(j) = d) = \Prob (d_{j}(i) = d)\Prob(d_n(j) = d)$.

Let $R_p(k)$ be the event that there is an edge from at least one vertex $p+i$ with $i \geq k$ to a vertex $p$. Then
\begin{multline*}
\Prob (d_n(i) = d, d_n(j) = d) - \Prob (d_n(i) = d)\Prob(d_n(j) = d)
\\
\le \Prob (d_j(i) = d, d_n(j) = d) + \Prob(R_i(j-i)) - \Prob (d_n(i) = d)\Prob(d_n(j) = d)
\\
= \Prob (d_j(i) = d)\Prob(d_n(j) = d) + \Prob(R_i(j-i)) - \Prob (d_n(i) = d)\Prob(d_n(j) = d)
\\
\le \Prob (d_j(i) = d)\Prob(d_n(j) = d)+ \Prob(R_i(j-i))\Prob(d_n(j) = d)
\\
+ \Prob(R_i(j-i)) - \Prob (d_n(i) = d)\Prob(d_n(j) = d)
\\
\le 2 \Prob(R_i(j-i)) = 2 \int_{a}^\infty R_i^q(j-i) f(q) dq
\\
= O\left( \int_{a}^{\infty} q^{-\gamma-1} e^{-\frac{j-i}{N}} dq \right)
= O\left( e^{-\frac{j-i}{N}} \right)\,.
\end{multline*}

Finally,
$$
Var[N_n(d)] = O\left( \sum_{1 \le i \le j \le n} e^{-\frac{j-i}{N}} \right) = O \left(N n\right)\,.
$$

Applying Chebyshev's inequality we get
$$
\Prob(|N_n(d) - \E [N_n(d)]| > \sqrt{Nn\log n}) = O\left(\frac{Var[N_n(d)]}{Nn\log n}\right) = O\left(\frac{1}{\log n}\right)\,.
$$

\subsection{Recency property}

In this section, we show that the behavior of $e(T)$ for the model with exponential popularity decay is realistic. It was shown in \cite{2012model} that $e(T)$ decreases exponentially with $T$ in real data.

First, we compute the expectation of $e(T)$. The following theorem holds.
\begin{thm}\label{th:recency2}
For any integer $T$
$$
\E [e(T)] = e^{-\frac{T}{N}} + O\left(\frac{N}{n}\right)\,.
$$
\end{thm}
Indeed, the probability that an edge from a vertex $k$ goes to a vertex $i$ with $k-i>T$ equals
$e^{-\frac{T}{N}} + O\left(e^{-\frac{k}{N}}\right)$. From this Theorem~\ref{th:recency2} follows.

As in the Section~\ref{sec:recency}, we can use Chebyshev's inequality to prove the concentration.

\begin{thm}\label{th:recency2conc}
For any integer $T$
$$
\Prob\left(|e(T) - \E [e(T)]| \ge  \sqrt{\frac{N \log n}{n}} \right) = O\left( \frac{1}{\log n} \right) \,.
$$
\end{thm}

\begin{Proof}

As before, by $e_i$ and $e_j$ we denote any two different edges drawn from vertices $i$ and $j$, and $l(e)$ is the age difference between the endpoints of an edge $e$.

Note that for $i \le j$
$$
\Prob(l(e_i) > T, l(e_j) > T|l(e_j)\le j-i) = \Prob(l(e_i) > T)\Prob(l(e_j) > T|l(e_j)\le j-i)\,.
$$
Therefore, we can estimate the following difference:
\begin{multline*}
\Prob(l(e_i) > T, l(e_j) > T) - \Prob(l(e_i) > T)\Prob(l(e_j) > T)
\\ = \Prob(l(e_j)\le j-i) \Big( \Prob(l(e_i) > T, l(e_j) > T| l(e_j)\le j-i)
\\ - \Prob(l(e_i) > T)\Prob(l(e_j) > T|l(e_j)\le j-i)\Big)
\\ + \Prob(l(e_j)> j-i) \Big( \Prob(l(e_i) > T, l(e_j) > T| l(e_j) > j-i)
\\ - \Prob(l(e_i) > T)\Prob(l(e_j) > T|l(e_j) > j-i)\Big) \le
\Prob(l(e_j)> j-i)
 \le e^{-\frac{j-i}{N}}\,.
\end{multline*}



Thus,
$$
Var[m\,n\,e(T)] = O\left( \sum_{1 \le i \le j \le n} e^{-\frac{j-i}{N}} \right) = O \left(N n\right)\,.
$$

Finally,
$$
\Prob\left(m\,n\,|e(T) - \E [e(T)]| \ge  m \sqrt{N \,n \log n} \right) \le \frac{Var[m\,n\,e(T)]}{m^2\, N\,n \log n} = O\left( \frac{1}{\log n} \right)\,.
$$

\end{Proof}

Theorems~\ref{th:recency2} and~\ref{th:recency2conc} mean that $e(T)$ decays exponentially, as it was observed in real data.

\section{Conclusion}

In this paper we analyze recency-based models.
The idea of adding the recency factor to the attractiveness function was first suggested in~\cite{2012model}.
In this paper we consider the most realistic model proposed in~\cite{2012model} and conduct a rigorous analysis of its properties.
In order to do this, we first provide a new formal definition of the model.
Then, we justify the fact that the power-law distribution of inherent quality leads to the power-law degree distribution.
We also rigorously analyze the recency property, i.e., the behavior of~$e(T)$, and prove that $e(T)$ decreases exponentially as it is observed in some real-world networks.

\section*{Funding}

This work was supported by the Russian Foundation for Basic Research [grant number 15-01-03530].

\section*{Acknowledgements}

The authors thank the anonymous reviewers for valuable feedback and suggestions.

\newpage

\section*{Appendix}

\subsection{Proof of Lemma~\ref{lem1}}

\textbf{Lemma~\ref{lem1}} \textit{
Let $\xi_1, \dots, \xi_n$ be mutually independent random variables, $\E [\xi_i] = 0$, $\E \left[|\xi_i|^{\alpha} \right] < \infty$, $1\le \alpha \le 2$, then
$$
\E \left[ |\xi_1 + \ldots + \xi_n|^{\alpha} \right] \le 2^{\alpha} \left(\E\left[|\xi_1|^{\alpha}\right]+\ldots+\E\left[|\xi_n|^{\alpha}\right] \right)\,.
$$}

\begin{Proof}

We use the following two facts.

\textbf{Fact 1.} If $\xi$ and $\eta$ are independent random variables and $\eta$ is symmetrically distributed, then for any $1 \le \alpha \le 2$
\begin{equation*}
\E \left[ |\xi + \eta|^\alpha \right] \le \E \left[|\xi|^\alpha\right] + \E \left[|\eta|^\alpha\right]\,.
\end{equation*}

\begin{Proof}
$$
\E \left[ |\xi + \eta|^\alpha \right] = \frac{1}{2} \left(\E\left[ |\xi+\eta|^\alpha \right] +\E\left[ |\xi-\eta|^\alpha\right] \right)\,
$$
and it remains to show that for any $x$, $y$, and $1 \le \alpha \le 2$ we have
$$
\frac{1}{2} \left(|x+y|^\alpha + |x - y|^\alpha \right) \le  |x|^\alpha + |y|^\alpha\,.
$$
Without loss of generality we assume that $x\ge y\ge 0$ and consider the function $f(x,y) = \frac{1}{2} \left((x+y)^\alpha + (x - y)^\alpha \right) - x^\alpha - y^\alpha$. In order to show that $f(x,y) \le 0$ we note that $f(x,0) = 0$ and $\frac{\partial f(x,y)}{\partial y} \le 0$. In turn, $\frac{\partial f(x,y)}{\partial y} \le 0$ since $\frac{\partial f(x,y)}{\partial y} \big|_{x=y} \le 0$ and $\frac{\partial^2 f(x,y)}{\partial x \, \partial y} \le 0$.
\end{Proof}

\textbf{Fact 2.}  If $\alpha \ge 1$, $\xi$ and $\eta$ are independent random variables, $\E \left[ \eta \right] = 0$, $\E \left[ |\xi|^\alpha \right]<\infty$, $\E \left[ |\eta|^\alpha \right]<\infty$, then
\begin{equation*}
\E \left[ |\xi + \eta|^\alpha \right] \ge \E \left[ |\xi|^\alpha \right]\,.
\end{equation*}

\begin{Proof}
Fact 2 follows directly from Jensen's inequality.
\end{Proof}

Now, let us prove Lemma~\ref{lem1}. Consider random variables $\xi_1', \dots, \xi_n'$, such that $\xi_i'$ has the same distribution as $\xi_i$ and $\xi_1, \dots, \xi_n, \xi_1', \dots, \xi_n'$ are mutually independent. Note that $\xi_i-\xi_i'$ is symmetrically distributed for any $i$. Then from Facts~1 and~2 it follows that
\begin{multline*}
\E \left[ |\xi_1 + \ldots + \xi_n|^{\alpha} \right] \le \E \left[ |\xi_1-\xi_1' + \ldots + \xi_n-\xi_n'|^{\alpha} \right] \\ \le \E \left[ |\xi_1-\xi_1'|^{\alpha} + \ldots + |\xi_n-\xi_n'|^{\alpha} \right]\,.
\end{multline*}
Finally, it remains to note that $\E\left[|\xi_i-\xi_i'|^{\alpha}\right] \le 2^{\alpha}\E\left[|\xi_i|^{\alpha}\right]$.

\end{Proof}

\subsection*{Proof of Lemma~\ref{S-+}}

First, recall the statement of Lemma~\ref{S-+} from Section~\ref{sec:expectation}.

\textbf{Lemma~\ref{S-+}} \textit{
Assume that both $N$ and $d$ grow but ${d = o(N^{1-c})}$. If $1/2 \le c \le 1$, then
$$
S_{\mp} =  \frac{ \gamma }{d^{\gamma + 1}} \left( \frac{(\gamma-1)m}{\gamma} \right)^\gamma
           \left( 1 + o(1) \right)\,.
$$}

\begin{Proof}

Recall that
$$
S_{\mp} = \int_{a}^{N^c/2} f(q) {mN \choose d} \left(\frac{q}{N\E[\zeta] \pm N^c}\right)^{d} \left(1 - \frac{q}{N\E[\zeta] \mp N^c}\right)^{mN-d} d q\,.
$$

Let us rewrite $S_{\mp}$ using the incomplete beta-function $B(x; a, b)$

\begin{multline*}
S_{\mp} = \int_{a}^{N^c/2}
\frac{\gamma a^\gamma}{q^{\gamma + 1} }
{mN \choose d} \left(\frac{q}{N\E[\zeta] \pm N^c}\right)^{d} \left(1 - \frac{q}{N\E[\zeta] \mp N^c}\right)^{mN-d} d q
\\
=   \frac{ \gamma a^\gamma {mN \choose d} \left(N\E[\zeta] \mp N^c\right)^{d-\gamma-1} }{\left(N\E[\zeta] \pm N^c\right)^{d}}  \\
\cdot \int_{a}^{N^c/2}
 \left(\frac{q}{N\E[\zeta] \mp N^c}\right)^{d - \gamma - 1} \left(1 - \frac{q}{N\E[\zeta] \mp N^c}\right)^{mN-d} d q \\
=   \frac{ \gamma a^\gamma {mN \choose d} \left(N\E[\zeta] \mp N^c\right)^{d-\gamma} }{\left(N\E[\zeta] \pm N^c\right)^{d}}
\cdot \int_{\frac{a}{N\E[\zeta] \mp N^c}}^{\frac{N^c/2}{N\E[\zeta] \mp N^c}}
 x^{d - \gamma - 1} \left(1 - x\right)^{mN-d} d x \\
=  \frac{ \gamma a^\gamma {mN \choose d} \left(N\E[\zeta] \mp N^c\right)^{d-\gamma} }{\left(N\E[\zeta] \pm N^c\right)^{d}}
 \left(
    B\left(\frac{N^c / 2}{N\E[\zeta] \mp N^c}; d - \gamma, mN-d+1\right)
     \right.\\ \left. -
    B\left(\frac{a}{N\E[\zeta] \mp N^c}; d - \gamma, mN-d+1\right)
 \right)\,. \\
\end{multline*}

Let us denote $\frac{N^c}{N\E[\zeta]}$ by $\varepsilon$, then we get



\begin{multline*}
S_{\mp} = \frac{ \gamma a^\gamma {mN \choose d}}{ \left(N\E[\zeta]\right)^{\gamma} }
           \frac{ \left(1 \mp \varepsilon \right)^{d-\gamma} }{\left(1 \pm \varepsilon\right)^{d}}
 \left(
    B\left(\frac{\varepsilon}{2 \mp 2 \varepsilon}; d - \gamma, mN-d+1\right)
     \right.\\ \left. -
    B\left(\frac{a}{ N\E[\zeta](1 \mp \varepsilon)}; d - \gamma, mN-d+1\right)
 \right)\,. \\
\end{multline*}

We will use the following estimates for the incomplete beta-function:

$$
B(x; a, b) = \int_0^x t^{a-1} (1-t)^{b-1} d t = O \left( \int_0^x t^{a-1} d t \right) = O \left( \frac{x^a}{a} \right)\,,
$$

$$
B(x; a, b) = B(a,b) -  \int_x^1 (1-t)^{b-1} d t  = B(a,b) + O \left( \frac{(1-x)^b}{b} \right)\,.
$$

These estimates give us
\begin{multline*}
S_{\mp} =  \frac{ \gamma a^\gamma {mN \choose d}}{ \left(N\E[\zeta]\right)^{\gamma} }
           \frac{ \left(1 \mp \varepsilon \right)^{d-\gamma} }{\left(1 \pm \varepsilon\right)^{d}}
 \left(
    B\left(d - \gamma, mN-d+1\right)  \right.\\ \left.
      +O \left( \frac{\left(1-\frac{\varepsilon}{2 \mp 2\varepsilon}\right)^{mN-d+1}}{mN-d+1} \right)
      +  O \left( \frac{\left(\frac{a}{N\E[\zeta] (1 \mp \varepsilon)}\right)^{d - \gamma}}{d - \gamma} \right)
 \right) \,.\\
\end{multline*}

We will use the fact that
${\frac{1}{B\left(d - \gamma, mN-d+1\right)}
 = \frac{\Gamma(mN+1-\gamma)}{\Gamma(d - \gamma)\Gamma(mN-d+1)}
 = O\left(\frac{(mN)^{d - \gamma}}{\Gamma(d - \gamma)}\right)}$
and factor out the beta-function:

\begin{multline}\label{eq:Smp}
S_{\mp} =  \frac{ \gamma a^\gamma {mN \choose d}}{ \left(N\E[\zeta]\right)^{\gamma} }
            B\left(d - \gamma, mN-d+1\right)
           \frac{ \left(1 \mp \varepsilon \right)^{d-\gamma} }{\left(1 \pm \varepsilon\right)^{d}}
   \\ \cdot
 \left(
       1
      +O \left( \frac{\left(1-\frac{\varepsilon}{2 \mp 2\varepsilon}\right)^{mN-d+1} (mN)^{d - \gamma}}{\Gamma(d - \gamma)(mN-d+1)} \right)
      +O \left( \frac{\left(\frac{a m}{\E[\zeta] (1 \mp \varepsilon)}\right)^{d - \gamma}}{\Gamma(d - \gamma + 1)} \right)
 \right) \,.\\
\end{multline}


Recall that $\varepsilon = N^{c-1}/E\zeta$. Let us simplify Equation~\eqref{eq:Smp}:

\begin{itemize}
\item[1.] $\frac{ \gamma a^\gamma {mN \choose d}}{ \left(N\E[\zeta]\right)^{\gamma} }
            B\left(d - \gamma, mN-d+1\right)
      = \frac{ \gamma a^\gamma }{\left(N\E[\zeta]\right)^{\gamma}}
        \frac{\Gamma(mN + 1)}{\Gamma(d + 1)\Gamma(mN-d+1)}
        \frac{\Gamma(d - \gamma)\Gamma(mN-d+1)}{\Gamma(mN+1-\gamma)}  $
      $ = \frac{ \gamma }{d^{\gamma + 1}} \left( \frac{am}{\E[\zeta]} \right)^\gamma (1 + o(1))$
if $d$ and $N$ grow.
\item[2.] $\frac{ \left(1 \mp \varepsilon \right)^{d-\gamma} }{\left(1 \pm \varepsilon\right)^{d}} =
       \frac{ \left(1 \mp N^{c-1}/E\zeta  \right)^{d-\gamma} }{\left(1 \pm N^{c-1}/E\zeta  \right)^{d}} = 1 + o(1)$ if $d = o(N^{1-c})$.
\item[3.] $
O \left( \frac{\left(1-\frac{\varepsilon}{2 \mp 2\varepsilon}\right)^{mN-d+1} (mN)^{d - \gamma}}{\Gamma(d - \gamma)(mN-d+1)} \right) =
O \left(
                      \frac{ e^{(mN-d+1)\log{\left(1-\frac{N^{c-1}}{2 \E[\zeta] \mp 2 N^{c-1}} \right)} + (d - \gamma - 1) \log(mN) } }
                           {\Gamma(d - \gamma)(1-\frac{d-1}{mN})}
        \right)
$

$
= O \left(
                      \frac{ e^{\frac{- m N^c \left(1-\frac{d-1}{mN} \right)}{2 E\zeta \mp 2 N^{c-1}} + (d - \gamma - 1) \log(mN) } }
                           {\Gamma(d - \gamma)(1-\frac{d-1}{mN})}
        \right)
= O \left(
                      \frac{ e^{\frac{- m  N^c \left(1-o(1)\right)}{2 E\zeta \mp o(1)} + d \log (mN) } }
                           {\Gamma(d - \gamma)(1-o(1))}
        \right)
= o(1) $ if $d < \frac{m N^c}{2  E\zeta \log (mN)}$ which is true for sufficiently large $N$ as soon as $d = o(N^{1-c})$ and ${1/2 < c < 1}$.

\item[4.]  $O \left( \frac{\left(\frac{a m}{\E[\zeta] (1 \mp \varepsilon)}\right)^{d - \gamma}}{\Gamma(d - \gamma + 1)} \right) = o(1)$
if $d$ and $N$ grow.
\end{itemize}

Thus, we first demonstrated that the main term in~\eqref{eq:Smp} is equal to $\frac{ \gamma }{d^{\gamma + 1}} \left( \frac{am}{\E[\zeta]} \right)^\gamma (1 + o(1))$. Then, we showed that the error multiplier equals $1+o(1)$. Finally, we proved that two error summands are equal to $o(1)$.
To conclude the proof of the lemma it remains to note that $\E[\zeta] = \frac{\gamma a}{\gamma-1}$, therefore $\frac{ \gamma }{d^{\gamma + 1}} \left( \frac{am}{\E[\zeta]} \right)^\gamma = \frac{ \gamma }{d^{\gamma + 1}} \left( \frac{(\gamma-1)m}{\gamma} \right)^\gamma$.

\end{Proof}

\subsection*{Proof of Lemma~\ref{rho-+_exp}}

First, recall the statement of Lemma~\ref{rho-+_exp} from Section~\ref{sec:expectation2}.

\textbf{Lemma~\ref{rho-+_exp}}\textit{
Assume that both $d$ and $N$ grow, $d = o(N^{1-c})$, $d = o(e^{\varphi(N)})$, and $q \leq N^r$, then
$$
S_{\mp}(d,q) = \frac{ \gamma }{d^{\gamma + 1}} \left( \frac{(\gamma-1)m}{\gamma} \right)^\gamma \left(1 + o(1)\right)\,.
$$}

\begin{Proof}

Recall that
$$
S_{\mp}(d,q) := \int_{a}^{N^r} f(q) \rho_{\mp}(d,q) d q\,,
$$
where
\begin{equation}\label{eq:rho1}
\rho_{\mp}(d,q) =
\prod_{i=0}^{m \varphi(N) N} \left( 1 - \frac{q e^{\frac{-[ i / m ]}{N}}}{N \E [\zeta] \pm N^c}  \right)
\sum_{0 \leq i_1 < \ldots < i_d \leq m \varphi(N) N}
\prod_{j=1}^{d} \frac{ \frac{q e^{\frac{-[ i_j / m ]}{N}}}{N \E [\zeta] \mp N^c}}
                     {1 - \frac{q e^{\frac{-[ i_j / m ]}{N}}}{N \E [\zeta] \pm N^c}} \,.
\end{equation}

Therefore, we first prove the following lemma on the behavior of $\rho_{\mp}(d,q)$.

\begin{lem}
Under the condition of Lemma~\ref{rho-+_exp} we have
$$
\rho_{\mp}(d,q) =
 \left(1 + O\left(\frac{q^2}{N}\right)+ o(1) \right)
 \left( \frac{q m}{\E [\zeta]} \right)^d  \frac{e^{\frac{-q m}{\E [\zeta]}}}{d!}\,.
$$
\end{lem}

\begin{Proof}


Note that
\begin{multline}\label{eq:part1}
\prod_{i=0}^{m \varphi(N) N} \left( 1 - \frac{q e^{\frac{-[ i / m ]}{N}}}{N \E [\zeta] \pm N^c}  \right) =
\prod_{i=0}^{m \varphi(N) N} \left( 1 - \frac{q e^{\frac{-i / m}{N}} e^{O\left(\frac{1}{N}\right)}}{N \E [\zeta] \pm N^c}  \right)
\\
= \exp \left( \sum_{i=0}^{m \varphi(N) N}
\log \left(
    1 - \frac{q e^{\frac{-i / m}{N}} \left(1 + O\left(\frac{1}{N}\right)\right)}{N \E [\zeta] \pm N^c}
\right) \right)
\\
= \exp \left( -\sum_{i=0}^{m \varphi(N) N}
\left(
   \frac{q e^{\frac{-i}{mN}} \left(1 + O\left(\frac{1}{N}\right)\right)}{N \E [\zeta] \pm N^c}
   + O \left( \frac{q^2 e^{\frac{-2i}{mN}}}{N^2} \right)
\right) \right)
\\
=
\exp \left( -\frac{q \left(1 - e^{-\varphi(N)}\right)}{(1 - e^{\frac{-1}{mN}})(N \E [\zeta] \pm N^c)} + \right. $$ $$ \left.
+ O \left( \frac{q \left(1 - e^{-\varphi(N)}\right)}{N^2 (1 - e^{\frac{-1}{mN}})} \right)
+ O \left( \frac{q^2 \left(1 - e^{-2\varphi(N)}\right)}{N^2 (1 - e^{\frac{-2}{mN}})} \right)
\right) \\
= \left(1 + O\left(\frac{q^2}{N}\right) + O\left(N^{c-1}\right) + O\left(e^{-\varphi(N)}\right)\right)
\exp \left(\frac{-q m}{\E [\zeta]} \right) \\
= \left(1 + o(1)\right)
\exp \left(\frac{-q m}{\E [\zeta]} \right) \,.
\end{multline}
Here we used the fact that $q \leq N^r < N^{1/2}$ since it allows us to
estimate $\exp \left(O\left( \frac{q^2}{N} \right) \right)$ as $1 + O \left( \frac{q^2}{N} \right)$.

Let us continue
\begin{multline}\label{eq:part2}
\prod_{j=1}^{d} \frac{q}{\left(1 - \frac{q e^{\frac{-[ i_j / m ]}{N}}}{N \E
[\zeta]\mp N^c}\right)\left(N\E[\zeta]\mp N^c\right)}
\\
=\left(\frac{q}{N\E[\zeta]}\right)^d \left(1 + O\left(\frac{d\,q}{N}\right)+O\left(dN^{c-1}\right)\right) =
\left(\frac{q}{N\E[\zeta]}\right)^d \left(1 + o(1)\right)\,,
\end{multline}
\begin{multline}\label{eq:part3}
\sum_{0 \leq i_1 < \ldots < i_d \leq m \varphi(N) N}
\prod_{j=1}^{d} e^{\frac{-[ i_j / m ]}{N}} =
\sum_{0 \leq i_1 < \ldots < i_d \leq m \varphi(N) N}
 e^{\frac{- i_1 - \ldots - i_d }{m N}} \left(1 + O\left(\frac d N\right) \right)
 \\
=
\sum_{0 \leq i_1 < \ldots < i_d \leq m \varphi(N) N}
 e^{\frac{- i_1 - \ldots - i_d }{m N}} \left(1 + o\left(1\right) \right)\,.
\end{multline}

It remains to estimate $\sum_{0 \leq i_1 < \ldots < i_d \leq m \varphi(N) N} e^{\frac{-i_1 - \ldots - i_d}{mN}}$.
We use the following notation:
$$
F(k,d) = \sum_{0 \leq i_1 < \ldots < i_d \leq m \varphi(N) N} e^{\frac{-i_1 - \ldots - i_{d-1} - k\,i_d}{mN}}\,.
$$

\begin{lem}\label{lem3}
If $d\,(k+d) = o(N)$ and $k+d = o\left(e^{\varphi(N)}\right)$, then
$$
F(k,d) = \frac{(mN)^d(k-1)!}{(k+d-1)!} \left(1 + o(1)\right)\,.
$$
\end{lem}

\begin{Proof}

Note that
$$
F(k,1) = \sum_{0 \leq i_1 \leq m \varphi(N) N} e^{\frac{-k\, i_1}{mN}} =
 \frac{1 - e^{-k\varphi(N)-\frac{k}{mN}}}{1 - e^{-\frac{k}{m N}}}\,.
$$
Let us get a recurrent formula for $F(k,d)$:
\begin{multline*}
F(k,d) = \sum_{0 \leq i_1 < \ldots < i_d \leq m \varphi(N) N} e^{\frac{-i_1 - \ldots - i_{d-1} - k\, i_d}{mN}}
\\
= \sum_{0 \leq i_1 < \ldots < i_{d-1} \leq m \varphi(N) N} e^{\frac{-i_1 - \ldots - i_{d-1}}{mN}}
\sum_{i_d = i_{d-1}+1}^{m \varphi(N) N} e^{\frac{-k\,i_d}{mN}}
\\
= \sum_{0 \leq i_1 < \ldots < i_{d-1} \leq m \varphi(N) N} e^{\frac{-i_1 - \ldots - i_{d-1}}{m N}}
   \frac{e^{-\frac{k(i_{d-1}+1)}{m N}} - e^{-k\varphi(N)-\frac{k}{mN}}}{1 - e^{-\frac{k}{m N}}}
   \\
=  \frac{e^{-\frac{k}{mN}}}{1 - e^{-\frac{k}{m N}}}
\sum_{0 \leq i_1 < \ldots < i_{d-1} \leq  N \varphi(N)} \left(  e^{\frac{-i_1 - \ldots - (k+1) i_{d-1}}{m N}} - e^{-k\varphi(N)} e^{\frac{-i_1 - \ldots - i_{d-1}}{m N}} \right)
\\
=  \frac{e^{-\frac{k}{mN}}}{1 - e^{-\frac{k}{m N}}}
 \left(  F(k+1,d-1) - e^{-k\varphi(N)} F(1,d-1) \right)\,.
\end{multline*}

It is easy to get an upper bound for $F(k,d)$
\begin{multline*}
F(k,d) \le \frac{e^{-\frac{k}{mN}}}{1 - e^{-\frac{k}{m N}}}
  F(k+1,d-1)  \le \ldots
  \\
\le \frac{e^{-\frac{(2k+d-2)(d-1)}{2mN}}F(k+d-1,1) }{\left(1 - e^{-\frac{k}{m N}}\right)\ldots\left(1 - e^{-\frac{k+d-2}{m N}}\right)}
   \le
\frac{e^{-\frac{(2k+d-2)(d-1)}{2mN}}}{\left(1 - e^{-\frac{k}{m N}}\right)\ldots\left(1 - e^{-\frac{k+d-1}{m N}}\right)}
\\
= \frac{e^{-\frac{(2k+d-2)(d-1)}{2mN}}}{\left(1 - e^{-\frac{k}{m N}}\right)\ldots\left(1 - e^{-\frac{k+d-1}{m N}}\right)}
= \frac{(mN)^d(k-1)!}{(k+d-1)!} \left( 1 + O\left(\frac{(k+d)d}{N}\right) \right)\,.
\end{multline*}

Using this upper bound and the recurrent formula above we can find a lower bound.
Assume that
\begin{multline*}
F(k,d) = \frac{e^{-\frac{k}{mN}}}{1 - e^{-\frac{k}{m N}}}
 \left(  F(k+1,d-1) - e^{-k\varphi(N)} F(1,d-1) \right) = \ldots
 \\
= \frac{e^{-\frac{(2k+d-2)(d-1)}{2mN}}F(k+d-1,1) }{\left(1 - e^{-\frac{k}{m N}}\right)\ldots\left(1 - e^{-\frac{k+d-2}{m N}}\right)} 
- \sum_{i=1}^{d-1} \frac{e^{-\frac{(2k+i-1)i}{2mN}} e^{-(k+i-1)\varphi(N)}F(1,d-i) }{\left(1 - e^{-\frac{k}{m N}}\right)\ldots\left(1 - e^{-\frac{k+i-1}{m N}}\right)}
\\
= \frac{(mN)^d(k-1)!}{(k+d-1)!} \left( 1 + O\left(\frac{(k+d)d}{N}\right) \right)
\\
- \sum_{i=1}^{d-1} e^{-(k+i-1)\varphi(N)} \frac{(mN)^{d-i}}{(d-i)!} \frac{(mN)^i (k-1)!}{(k+i-1)!} \\ \cdot \left(1 + O\left( \frac{(d-i)^2}{N}\right) + O\left(\frac{(k+i)i}{N} \right) \right)
  \\
= \frac{(mN)^d (k-1)!}{(k+d-1)!} \left(  1 + o(1)
- \sum_{i=1}^{d-1} \frac{e^{-(k+i-1)\varphi(N)} (k+d-1)! }{(d-i)! (k+i-1)!} \left(1 + o(1)\right)
\right)
  \\
\ge \frac{(mN)^d (k-1)!}{(k+d-1)!} \left(  1 + o(1) -
 \left(1 + o(1) \right)\sum_{i=1}^{d-1} \frac{ \left( (k+d-1) e^{-\varphi(N)}  \right)^{k+i-1} }{(k+i-1)!}  \right)
 \\
= \frac{(mN)^d (k-1)!}{(k+d-1)!} \left(  1 + o(1) +
O \left(\frac{ \left( (k+d-1) e^{-\varphi(N)}  \right)^{k} }{k!}  \right)
\right)
\\
= \frac{(mN)^d (k-1)!}{(k+d-1)!} \left(  1 + o(1)
\right) \,.
\end{multline*}

\end{Proof}


Finally, taking into account Equations~\eqref{eq:rho1}-\eqref{eq:part3} and Lemma~\ref{lem3}, we get
\begin{multline*}
\rho_{\mp}(d,q) = e^{\frac{-q m}{\E [\zeta]}} \left(\frac{q}{N\E[\zeta]}\right)^d F(1,d) (1+o(1))
\\
 =e^{\frac{-q m}{\E [\zeta]}} \frac{(mN)^d}{d!} \left(\frac{q}{N\E[\zeta]}\right)^d (1+o(1)) =
 \left(1 + o(1) \right)
 \left( \frac{q m}{\E [\zeta]} \right)^d  \frac{e^{\frac{-q m}{\E [\zeta]}}}{d!}\,.
\end{multline*}

\end{Proof}

Now we can estimate
$$
S_{\mp}(d,q) = \int_{a}^{N^r} f(q) \rho_{\mp}(d,q) d q 
=  \int_{a}^{N^r}
\left(1  + o(1)\right)
\frac{\gamma a^{\gamma}}{q^{\gamma+1}}
\left( \frac{q m}{\E [\zeta]} \right)^d  \frac{e^{\frac{-qm}{\E [\zeta]}}}{d!} d q \,.
$$

We get an incomplete gamma function:
\begin{multline*}
S_{\mp}(d,q) =
\left(1  + o(1)\right) \frac{\gamma a^{\gamma} m^\gamma}{d! ( \E [\zeta] )^{\gamma}}
\int_{\frac{am}{\E [\zeta]}}^{\frac{N^r m}{\E [\zeta]}}
x^{d-\gamma-1} e^{-x} d x
\\
=  \frac{\gamma(1 + o(1))}{\Gamma(d+1)}  \left( \frac{am}{\E [\zeta]} \right)^{\gamma}
\left(
      \Gamma\left(d-\gamma, \frac{am}{\E [\zeta]}\right)
       - \Gamma\left(d-\gamma, \frac{N^r m}{\E [\zeta]}\right) \right)
\\
=   \left(1 + o(1) \right)  \frac{\gamma \Gamma(d-\gamma)}{\Gamma(d+1)} \left( \frac{am}{\E [\zeta]} \right)^{\gamma}
=   \frac{\gamma}{d^{\gamma+1}} \left( \frac{am}{\E [\zeta]} \right)^{\gamma} \left(1 + o(1)\right)
\\
= \frac{ \gamma }{d^{\gamma + 1}} \left( \frac{(\gamma-1)m}{\gamma} \right)^\gamma(1+o(1))\,.
\end{multline*}

This concludes the proof.

\end{Proof}

\newpage

\bibliography{references}
\bibliographystyle{abbrv}

\label{lastpage}

\end{document}